\documentclass[amsppt,11pt]{amsart}
\usepackage{curves}
\usepackage{amsmath}

\newtheorem {theorem}{Theorem}[section]

\newtheorem {lemma}[theorem]{Lemma}

\newcounter{conjecture}
\newcounter{remark}

\newenvironment{remark}{\medskip{\bf Remark
\theremark.}
\addtocounter{remark}{1}}{}
\newcommand{\eqnsection}{
  \renewcommand{\theequation}{\thesection.\arabic{equation}}
   \makeatletter
   \csname @addtoreset\endcsname{equation}{section}
   \makeatother}

\newcommand{\be}{{\begin{equation}}}
\newcommand{\ee}{{\end{equation}}}
\def \bt{\begin{theorem}}
\def \et{\end{theorem}}
\def \bea{\begin{eqnarray}}
\def \eea{\end{eqnarray}}
\def \bas{\begin{eqnarray*}}
\def \eas{\end{eqnarray*}}

\def \al{\alpha}
\def \bb{\beta}
\def \ga{\gamma}
\def \Ga{\Gamma}
\def \de{\delta}

\def \ep{\epsilon}

\def \la{\lambda}

\def \om{\omega}
\def \Om{\Omega}

\def \si{\sigma}

\def \th{\theta}

\def \ff{\infty}

\def \wt{\widetilde}

\def \rar{\rightarrow}

\newcommand{\ls}[1]
   {\dimen0=\fontdimen6\the\font \lineskip=#1\dimen0
\advance\lineskip.5\fontdimen5\the\font
\advance\lineskip-\dimen0
\lineskiplimit=.9\lineskip \baselineskip=\lineskip
\advance\baselineskip\dimen0 \normallineskip\lineskip
\normallineskiplimit\lineskiplimit
\normalbaselineskip\baselineskip
\ignorespaces }
\newcommand{\req}[1]{(\ref{#1})}

\def \Z{{\Bbb{Z}}}

\def \FF{{\mathcal F}}
\def \GG{{\mathcal G}}
\def \HH{{\mathcal H}}

\def \({\left(}
\def \){\right)}

\def \lc{\left\{}
\def \rc{\right\}}

\def \nn{\nonumber}

\def \bc{\begin{center} }
\def \ec{\end{center} }
\begin{document}

\eqnsection
\newcommand{\Ini}{{I_{n,i}}}
\newcommand{\reals}{{\Bbb{R}}}
\newcommand{\F}{{\mathcal F}}
\newcommand{\D}{{\mathcal D}}
\newcommand{\Fn}{{{\mathcal F}_n}}
\newcommand{\Gn}{{{\mathcal G}_n}}
\newcommand{\Hn}{{{\mathcal H}_n}}
\newcommand{\Fp}{{{\mathcal F}^p}}
\newcommand{\Gp}{{{\mathcal G}^p}}
\newcommand{\PPP}{{P}}
\newcommand{\Pop}{{P\otimes \PPP}}
\newcommand{\hm}{\HH^\varphi}
\newcommand{\nuw}{{\nu^W}}
\newcommand{\ths}{{\theta^*}}
\newcommand{\beq}[1]{\begin{equation}\label{#1}}
\newcommand{\eeq}{\end{equation}}
\newcommand{\integers}{{\rm I\!N}}
\newcommand{\E}{{E}}
\newcommand{\te}{{\tilde{\delta}}}
\newcommand{\tI}{{\tilde{I}}}
\newcommand{\loge}{{\log(1/\ep)}}
\newcommand{\logen}{{\log(1/\ep_n)}}
\newcommand{\epn}{{\ep_n}}
\def\var{{\rm Var}}
\def\cov{{\rm Cov}}
\def\one{{\bf 1}}
\def\leb{{\mathcal L}eb}
\def\Ho{{\mbox{\sf H\"older}}}  
\def\thi{{\mbox{\sf Thick}}}
\def\cthi{{\mbox{\sf CThick}}}
\def\late{{\mbox{\sf Late}}}
\def\clate{{\mbox{\sf CLate}}}
\newcommand{\ffrac}[2]
  {\left( \frac{#1}{#2} \right)}
\newcommand{\calF}{{\mathcal F}}
\newcommand{\dfn}{\stackrel{\triangle}{=}}
\newcommand{\beqn}[1]{\begin{eqnarray}\label{#1}}
\newcommand{\eeqn}{\end{eqnarray}}
\newcommand{\oo}{\overline}
\newcommand{\uu}{\underline}
\newcommand{\bfcdot}{{\mbox{\boldmath$\cdot$}}}
\newcommand{\Var}{{\rm \,Var\,}}
\def\squarebox#1{\hbox to #1{\hfill\vbox to #1{\vfill}}}
\renewcommand{\qed}{\hspace*{\fill}
           
\vbox{\hrule\hbox{\vrule\squarebox{.667em}\vrule}\hrule}\smallskip}
\newcommand{\half}{\frac{1}{2}\:}
\newcommand{\beaa}{\begin{eqnarray*}}
\newcommand{\eeaa}{\end{eqnarray*}}
\newcommand{\calK}{{\mathcal K}}
\def\dimm{{\overline{{\rm dim}}_{_{\rm M}}}}
\def\dimp{\dim_{_{\rm P}}}
\def\htaum{{\hat\tau}_m}
\def\htaumk{{\hat\tau}_{m,k}}
\def\htaumkj{{\hat\tau}_{m,k,j}}
\def \IJK{\mathcal I}
\def \bga{\bar\gamma}
\def\taup{{\oo{\tau}}}
\def\sip{{\oo{\zeta}}}
\def\sii{{\zeta}}

\bibliographystyle{amsplain}

\title[Proof of the Erd\H{o}s-Taylor Conjecture] {A Random Walk
Proof of the Erd\H{o}s-Taylor Conjecture}

\author[ Jay Rosen] {Jay Rosen$^\ddagger$}

\date{March 2, 2005.
\newline\indent
$^\ddagger$Research supported, in part, by grants from the NSF  and from
PSC-CUNY }

\begin{abstract}
\noindent For the simple random walk in
$\Z^2$ we study those points which are
 visited an unusually large number of times, and provide a new proof of the 
Erd\H{o}s-Taylor Conjecture describing the number of visits to the most visited
point. 
\end{abstract}

\maketitle

\section{Introduction}

In our paper \cite{DPRZ4} we proved a conjecture of Erd\H{o}s and Taylor
concerning the number   $L^{ \ast}_n$ of visits to the most visited site for the
simple random walk in
$\Z^2$ up to step $n$. 
\begin{theorem}\label{theo-et}
\label{theo-late} Let $\{ X_{ j}\,;\,j\geq 1\}$ be the simple random walk in $\Z^2$.
 Then 
\begin{equation}
\lim_{n\to \ff}{L^{ \ast}_n\over (\log n)^2}=1/\pi
\hspace{.1in}\mbox{a.s.}\label{0.0h}
\end{equation}
\end{theorem}
 Our approach in that paper was to first prove an analogous result for planar 
Brownian motion and then to use strong approximation. The goal of
this paper, which is purely expository, is to show how
 to prove (\ref{0.0h}) using only random walk methods.  We also go
beyond (\ref{0.0h}) and study the size of the set of `frequent points', i.e. those
points in 
$\Z^2$ which are
 visited an unusually large number of times, 
 of order $(\log n)^2$. Perhaps more important than our specific results, we develop
powerful estimates for the simple random walk which we expect will have wide
applicability.

  Let $L_n^{ x}$ denote the number of times that $x\in \Z^2$ is visited by
 the random walk in $\Z^2$ up to step $n$, (so that  $L^{ \ast}_n=\max_{ x\in
\Z^2}L_n^{ x}$).  Set
\begin{equation}
\Theta_n(\al)=\Big\{x\in  \Z^2:\;
 \frac{L^{ x}_{ n}}{(\log n)^2}\geq \al/\pi
\Big\}.\label{bthick.1}
\end{equation}

Let $D(x,r)=\{y\in\Z^2\,|\,|y-x|<r\} $.  Let
$T_{B}=\inf\{i\geq 0\,|\,X_i\in B\}$ for any
$B\subseteq \Z^2$ and set
\begin{equation}
\Psi_n(a)=\Big\{x\in  D(0, n):\;
 \frac{L^{ x}_{ T_{D(0, n)^{ c}}}}{(\log n)^2}\geq 2a/\pi
\Big\}.\label{bthick.1}
\end{equation}

\begin{theorem}\label{theo-mf}
\label{theo-late} Let $\{ X_{ j}\,;\,j\geq 1\}$ be the simple random walk in $\Z^2$.
Then for any
$0<\al< 1$
\begin{equation}
\lim_{n\to \ff}{\log |\Theta_n(\al)|\over \log n}=1-\al
\hspace{.1in}\mbox{a.s.}\label{thick.2}
\end{equation} Equivalently, for any $0< a< 2$
\begin{equation}
\lim_{n\to \ff}{\log |\Psi_n(a)|\over \log n}=2-a
\hspace{.1in}\mbox{a.s.}\label{bthick.2}
\end{equation}
\end{theorem}

The equivalence of (\ref{thick.2}) and (\ref{bthick.2}) follows from the strong
invariance principle which implies that 
\begin{equation}
 \lim_{n\to \ff}\frac{\log T_{D(0, n)^{ c}}}{\log n}=2
\hspace{.1in}\mbox{a.s.}\label{equiv.1}
\end{equation}

In particular, (\ref{0.0h}) is equivalent to
\begin{equation}
\lim_{n\to \ff}{L^{ \ast}_{D(0, n)^{ c}}\over (\log n)^2}=2/\pi
\hspace{.1in}\mbox{a.s.}\label{0.0hk}
\end{equation}

In Section \ref{sec-upperbound} we prove the upper bound for Theorem
\ref{theo-mf}, and the lower bound is proven in Section \ref{sec-lowerbound},
subject to Lemma \ref{pmomlb}. Sections \ref{51momest}-\ref{harn} are devoted to
the proof of Lemma \ref{pmomlb}.

Here are the basic ideas behind our proof of (\ref{0.0hk}). It is not too hard to get
good estimates on $P( L^{ x}_{ T_{D(0, n)^{ c}}}\geq (2-\ep)(\log n)^2/\pi
)$ for $x\in D(0, n)$, and indeed this is how we obtain the upper bound for 
(\ref{0.0hk}) in Section \ref{sec-upperbound}. On the other hand, the lower bound
 requires a second moment estimate,  which is problematic due to the large
correlation between the events $\{ L^{ x}_{ T_{D(0, n)^{ c}}}\geq (2-\ep)(\log
n)^2/\pi\}$ for different $x\in D(0, n)$. Rather than study the events $\{ L^{ x}_{
T_{D(0, n)^{ c}}}\geq (2-\ep)(\log n)^2/\pi\}$, we define a point $x\in D(0, n)$ to
be `n-succesful' if, up to time $T_{D(0, n)^{ c}}$, there is a specified number of
excursions between circles centered  at $x$. The `excursion count' is chosen to be
`typical' for $x$ with $\{ L^{ x}_{
T_{D(0, n)^{ c}}}\geq (2-\ep)(\log n)^2/\pi\}$, so that  the
probability of $x$ being `n-succesful' is close to 
$P( L^{ x}_{ T_{D(0, n)^{ c}}}\geq (2-\ep)(\log n)^2/\pi)$ and also, for large
enough $n$, all `n-succesful' points $x$ have $L^{ x}_{ T_{D(0, n)^{ c}}}\geq
(2-2\ep)(\log n)^2/\pi$, see Lemma \ref{plem1}. The advantage of working with 
`n-succesful' points is that the correlation between two such points, $x,y$, can be
controlled using the tree-like structure of circles centered at the two points. In
particular we show that for $x,y$ far apart, excursions between circles centered at
$x$ and close to $x$ are almost independent of the excursions between circles
centered at
$y$ and close to $y$.

\section{Local time estimates and upper bounds}
\label{sec-upperbound}

Let $X_n,\,n\geq 0$ denote the simple random walk in $\Z^2$. We use 
$P^x,\,E^x$ for the probability and expectation for the walk started at $x$, and
write simply $P,\,E$ for the walk started at the
origin.  For any set
$A\subseteq
\Z^2$ we define the boundary $\partial A$ of $A$ by $\partial
A=\{y\in\Z^2\,|\,y\in A^c,\,\mbox{and }
\inf_{x\in A}|y-x|\leq 1\} $.  For
$x,y\in \Z^2$ define the Green's function
\begin{equation} G_A(x,y)=\sum_{i=0}^\ff P^x\(X_i=y,\,i<T_{A^c}\).\label{p1.1}
\end{equation}

\begin{lemma}\label{lem-hit}  For $|x_0|=r<R$ and some finite constant $\ga$,
\begin{equation}
\label{tohelet}
\E^{x_0}(L^{ 0}_{T_{D(0, R)^{ c}}})=
 G_{D(0,R)}(x_0,0)=\left\{ \begin{array}{ ll}{2\over \pi}\log (R) +\ga+O(
R^{-1}) &
\mbox{  if
$x_0=0$}\\ \\ {2\over \pi}\log ({R\over r})+O( r^{-1}) & \mbox{  if
$x_0\neq 0$.}
\end{array}
\right.
\end{equation} 
Let $x_0\neq 0$. For any $0<\varphi\leq 1$
\begin{equation}
E^{x_0}\(e^{ -{\varphi \over G_{D(0,R) }(0,0)}L^{ 0}_{T_{D(0, R)^{
c}}}}\)=1-{\log ({R\over r}) \over \log (R)}{\varphi\over 1+\varphi
}\(1+O({ 1\over \log (r)} )\)\label{babe.1}
\end{equation}
and for all
$z\geq 1$
\begin{equation}
 P^{x_0}( L^{ 0}_{T_{D(0, R)^{ c}}}\geq zG_{D(0,R)}(0,0))\leq c\sqrt{
z}e^{ -z}\,\label{prob}
\end{equation}
for some $c<\ff$ independent of $x_0, z, R$.
\end{lemma}

\noindent {\bf Proof of Lemma~\ref{lem-hit}:} By Theorem 1.6.6 of \cite{L1},
\begin{equation} G_{D(0,R) }(0,0)={ 2\over \pi}\log R +\ga+O(R^{ -1})
\label{p2.5}
\end{equation}
for an explicit constant $\ga$, and by Proposition 1.6.7 of \cite{L1}
\begin{equation} G_{D(0,R)}(x,0)={2\over \pi}\log \({R\over |x|}\)+O(
|x|^{-1}).
\label{p1.2}
\end{equation} 
Since
\begin{equation} L^{ 0}_{T_{D(0, R)^{ c}}}=\sum_{i<T_{D(0, R)^{ c}}} 1_{\{X_i=0\}
},\label{c1.1}
\end{equation} (\ref{tohelet}) follows.

For (\ref{babe.1}), note that, conditional  on  hitting $0$,  $L^{ 0}_{T_{D(0, R)^{
c}}}$ is a geometric random variable with  mean $G_{D(0,R) }(0,0)$.   Hence, 
\bea
&&
 E^{x_0}\(e^{ -{\varphi \over G_{D(0,R) }(0,0)}L^{ 0}_{T_{D(0,
R)^{ c}}}}\)
\nn
\\
&& =1-P^{ x_0}\(T_{ 0}<T_{D(0, R)^{c}}\)\nn\\
&&\hspace{ .5in}
 +P^{ x_0}\(T_{ 0}<T_{D(0, R)^{c}}\)\({ 1\over
 (e^{ {\varphi \over G_{D(0,R) }(0,0)}}-1 )G_{D(0,R) }(0,0)+1 }
\).\label{8.15}
\eea

Since by (\ref{p2.5})
\begin{equation} { 1\over G_{D(0,R) }(0,0)}=O(1/
\log (R))\label{8.15a}
\end{equation} 
we have 
\begin{equation}
(e^{ {\varphi \over G_{D(0,R) }(0,0)}}-1 )G_{D(0,R) }(0,0)+1
=1+\varphi +O({ 1\over \log (R)} ).\label{8.15aa}
\end{equation}

Furthermore, using the strong Markov property at the stopping time $T_{ 0}$
\begin{eqnarray}
G_{D(0,R) }(x_0,0)&=&P^{ x_0}\(\sum_{i<T_{D(0, R)^{ c}}} 1_{\{X_i=0\}
}\)\nn\\ 
&=&P^{ x_0}\(\sum_{i<T_{D(0, R)^{ c}}} 1_{\{X_i=0\}}\circ\theta_{ T_{ 0}};\,
T_{ 0}<T_{D(0, R)^{c}}\)\label{8.15f}\\ 
&=&P^{ x_0}\(T_{ 0}<T_{D(0, R)^{c}}\)P^{ 0}\(\sum_{i<T_{D(0, R)^{ c}}} 1_{\{X_i=0\}
}\)
\nonumber\\ 
&=&P^{ x_0}\(T_{ 0}<T_{D(0, R)^{c}}\)G_{D(0,R) }(0,0)
\nonumber
\end{eqnarray}
so that
\begin{equation}
P^{ x_0}\(T_{ 0}<T_{D(0, R)^{c}}\)={G_{D(0,R) }(x_0,0) \over G_{D(0,R)
}(0,0)}.\label{8.15fc}
\end{equation}

By (\ref{p2.5}) and (\ref{p1.2}),
\begin{equation} { G_{D(0,R) }(x_0,0)\over G_{D(0,R) }(0,0)}
={\log ({R\over r}) \over \log (R)}\(1+O({ 1\over \log (r)} )\),\label{8.15b}
\end{equation} 
and (\ref{babe.1}) then follows.

For (\ref{prob}) note that
by the strong Markov property   we have
\begin{eqnarray}
\E^{x_0}(L^{ 0}_{T_{D(0, R)^{ c}}})^k & \!\!\leq \!\!& k!
\E^{x_0}\(\sum_{0 \leq j_1 \leq \cdots\leq j_k < T_{D(0, R)^{ c}}}
\prod_{i=1}^k 1_{\{X_{j_{ i} }=0\} } \) \nn \\ &\!\!\leq \!\!& k!
\E^{x_0}\(\sum_{0 \leq j_1 \leq \cdots\leq j_{ k-1} <  T_{D(0, R)^{ c}}}
\prod_{i=1}^{ k-1} 1_{\{X_{j_{ i} }=0\} } G_{D(0,R)}(0,0)\)
\nn \\ &\!\!=\!\!&  k!G_{D(0,R)}(0,0)
\E^{x_0}\(\sum_{0 \leq j_1 \leq \cdots\leq j_{ k-1} <  T_{D(0, R)^{ c}}}
\prod_{i=1}^{ k-1} 1_{\{X_{j_{ i} }=0\} } \),
\nn
\end{eqnarray} 
and by iteration we obtain
\begin{equation}
\label{secmom}
\E^{x_0}(L^{ 0}_{T_{D(0, R)^{ c}}})^k\leq k!G_{D(0,R)}(x_0,0)(  G_{D(0,R)}(0,0))^{k
-1}.
\end{equation}
 To prove (\ref{prob}), use (\ref{secmom}), the fact that
$G_{D(0,R)}(x_0,0)\leq   G_{D(0,R)}(0,0)$  by (\ref{8.15f}), and Chebysheff to
obtain
\begin{equation}
\label{prob1} P^{x_0}( L^{ 0}_{T_{D(0, R)^{ c}}}\geq zG_{D(0,R)}(0,0))\leq { k!\over
z^{ k}}\,
\end{equation} then take $k=[z]$ and use Stirling's formula.
\qed

We next provide the required upper bounds in Theorem
\ref{theo-late}. Namely,  we will show that for any $a \in (0,2]$
\begin{equation}
\label{3.1b}
\limsup_{m\rar \ff}{ \log \Big|\Big\{x\in  D(0, m):\; L^{ x}_{ T_{D(0, m)^{ c}}}\geq
(2a/\pi) (\log m)^2\Big\}\Big|\over \log m}\leq 2-a \hspace{.2in}a.s.
\end{equation}

To see this fix $\de>0$ and note that by (\ref{prob}), for some $0<\ep<\de$, all
$x\in D(0, m)$ and all large enough $m$
\begin{equation}
\PPP\(  \frac{L^{ x}_{ T_{D(x, 2m)^{ c}}}}{(\log m)^2}\geq 2a/\pi \)\leq
m^{ -a+\ep}\label{prob2}
\end{equation}
 Therefore
\begin{eqnarray} &&
\PPP\(\Big|\Big\{x\in  D(0, m):\;
 \frac{L^{ x}_{ T_{D(0, m)^{ c}}}}{(\log m)^2}\geq 2a/\pi \Big\}\Big|\geq
m^{2-a+\de}\)\label{3.1ba}\\ &&\leq m^{-(2-a)-\de}\E^{0}\(\Big| \Big\{x\in  D(0,
m):\;
 \frac{L^{ x}_{ T_{D(0, m)^{ c}}}}{(\log m)^2}\geq 2a/\pi\Big\}\Big|\)\nn\\
&&=m^{-(2-a)-\de}\sum_{x\in  D(0, m)}\PPP^{0}\(  \frac{L^{ x}_{ T_{D(0, m)^{
c}}}}{(\log m)^2}\geq 2a/\pi
\)\nn\\ &&\leq m^{-(2-a)-\de}\sum_{x\in  D(0, m)}\PPP^{0}\( 
\frac{L^{ x}_{ T_{D(x, 2m)^{ c}}}}{(\log m)^2}\geq 2a/\pi \)\nn\\
 &&\leq m^{-(\de-\ep )}.\nn
\end{eqnarray} 
Now apply our result to $m=m_{ n}=e^{ n}$ to see by Borel-Cantelli
that for some $N(\om)<\ff$ a.s. we have that for all $n\geq N(\om)$
\begin{equation}
 \Big| \Big\{x\in  D(0, e^{ n}):\;\frac{L^{ x}_{ T_{D(0, e^{ n})^{ c}}}}{(\log e^{
n})^2}\geq 2a/\pi \Big\}\Big|\leq e^{(2-a+\de)n}.\label{prob3}
\end{equation}
Then if $e^{ n}\leq m\leq e^{ n+1}$
\begin{eqnarray} &&
 \Big| \Big\{x\in  D(0, m):\;\frac{L^{ x}_{ T_{D(0, m)^{ c}}}}{(\log m)^2}\geq
2a/\pi
\Big\}\Big|\label{prob4}\\ &&
 \leq \Big| \Big\{x\in  D(0, e^{ n+1}):\;\frac{L^{ x}_{ T_{D(0, e^{ n+1})^{ c}}}}{(\log
e^{ n})^2}\geq 2a/\pi
\Big\}\Big|\nn\\ && = \Big| \Big\{x\in  D(0, e^{ n+1}):\;\frac{L^{ x}_{ T_{D(0, e^{
n+1})^{ c}}}}{(\log e^{ n+1})^2}\geq 2a( 1+1/n)^{- 2}/\pi \Big\}\Big|\nn\\
&&
 \leq e^{(2-a( 1+1/n)^{ -2}+\de)(n+1)}\leq m^{(2-a( 1+1/n)^{
-2}+\de)(n+1)/n}.\nn
\end{eqnarray} 
(\ref{3.1b}) now follows on taking $\de\rar 0$.
\qed

\section{Lower bounds for probabilities}\label{sec-lowerbound}

Fixing $0<a<2$, we prove in this section that
\begin{equation}
\label{p.1}
\liminf_{m\rar \ff}{ \log \Big|\Big\{x\in  D(0, m):\; L^{ x}_{ T_{D(0, m)^{ c}}}\geq
(2a/\pi) (\log m)^2\Big\}\Big|\over \log m}\geq 2-a \hspace{.2in}a.s.
\end{equation} In view of (\ref{3.1b}), we will obtain Theorem \ref{theo-late} .

Set
$K_n=16e^{ n}n^{ 3 n}$. Using the fact that $\lim_{ n\rar
\ff}\log K_n/\log K_{ n-1}=1$, a simple interpolation argument shows that
 in order to prove (\ref{p.1}) it suffices to show that
\begin{equation}
\label{q.1}
\liminf_{n\rar \ff}{ \log \Big|\Big\{x\in  D(0, K_n):\; L^{ x}_{ T_{D(0, K_n)^{
c}}}\geq (2a/\pi) (\log K_n)^2\Big\}\Big|\over \log K_n}\geq 2-a
\hspace{.2in}a.s.
\end{equation} It suffices to prove that for any
$\de>0$ and sufficiently large $n$
\begin{equation}
\PPP\(\Big|\Big\{x\in  D(0, K_n):\; {L^{ x}_{ T_{D(0, K_n)^{c}}} \over  (\log K_n)^2}
\geq (2a/\pi)\Big\}\Big| \geq K_n^{2-a-\de }\)
\geq p_{ \de}>0.\label{q.2}
\end{equation} 
For then
\begin{equation}
\PPP\(\Big|\Big\{x\in  D(0, K_n):\; {L^{ x}_{ T_{D(0, K_n)^{c}}} \over   (\log K_{
n+1})^2}
\geq (2a-\de)/\pi\Big\}\Big| \leq K_{ n+1}^{2-a-2\de }\)
\leq1-p_{\de},\label{q.2ref}
\end{equation}
By considering the stopping times
$ \tau_{ l}=:T_{D(0, lK_n)^{c}}\,;\,l=0,1,2,\ldots, n^{ 3}-1$ 
and setting
\begin{equation}
A_{ l}=\lc\Big|\Big\{x\in  D(X_{  \tau_{ l}}, K_n):\; {L^{ x}_{ T_{D(X_{  \tau_{ l}},
K_n)^{c}}} \over   (\log K_{ n+1})^2}
\geq (2a-\de)/\pi\Big\}\Big| \leq K_{ n+1}^{2-a-2\de }\rc\label{q.2fr}
\end{equation}
we have that
\begin{eqnarray} &&
\lc\Big|\Big\{x\in  D(0, K_{ n+1}):\; {L^{ x}_{ T_{D(0, K_{ n+1})^{c}}} \over  
(\log K_{ n+1})^2}
\geq (2a-\de)/\pi\Big\}\Big| \leq K_{ n+1}^{2-a-2\de }\rc\label{q.2f}\\ &&
\hspace{ 2.5in}\subseteq \cap_{ l=0}^{ n^{ 3}-1}\,\,A_{ l}
\circ\th_{  \tau_{
l}}.\nonumber
\end{eqnarray}
Hence using the strong Markov property and (\ref{q.2ref}) we see that
\[
\PPP\(\Big|\Big\{x\in  D(0, K_{ n+1}):\; {L^{ x}_{ T_{D(0, K_{ n+1})^{c}}} \over  
(\log K_{ n+1})^2}
\geq (2a-\de)/\pi\Big\}\Big| \leq K_{ n+1}^{2-a-2\de }\)
\leq (1-p_{ \de})^{n^{ 3} }.\]
 An application of the Borel-Cantelli lemma followed by taking the
$\de\rar 0$ limit then gives (\ref{q.1}).

We start by constructing a subset of the set appearing in (\ref{q.2}), the probability
of which is easier to bound below. To this end fix $n$, and let 
$r_{ n,k}=e^{ n} n^{ 3( n-k)},\,k=0,\ldots, n$. In particular, $r_{ n,n}=e^{ n}$ and
$r_{ n,0}=e^{ n}n^{ 3n}$. Set
$K_n=16r_{ n,0}= 16e^{ n}n^{ 3n}$.   

Let $U_{ n}=[2r_{ n,0},3r_{ n,0}]^{ 2}\subseteq D(0, K_n)$. 
 For $x
\in U_{ n}$, let $N_{n,k}^x$ denote the number of excursions from
$\partial D(x,r_{n,k-1})$ to
$\partial D(x,r_{n,k})$ until time $T_{D(0,K_n)^{ c}}$. Set
$\frak{n}_k=3ak^2\log k$, and $k_{ 0}=4\vee\inf\{ k\,|\,\frak{n}_k\geq 2k\}$.

 We will say that a point
$x\in U_{ n}$ is {\bf $n$-successful} if 
 $N^x_{n,k}=1\,,\,\forall k=1,\ldots,k_{ 0}-1$ and 
\begin{equation} \frak{n}_k-k\leq N^x_{n,k}\leq \frak{n}_k+k\hspace{.3in}\forall
k=k_{ 0},\ldots,n
\label{pmperf}
\end{equation}

 Let
$Y(n,x)\,;\,x\in U_{ n}$ be the collection of random variables defined by
\[Y(n,x)=1\hspace{.1in}\mbox{if $x$ is $n$-successful}\] and $Y(n,x)=0$ otherwise.
Set
$\bar{q}_{n,x}=\PPP(Y(n,x)=1)=\E(Y(n,x))$.

 The next lemma relates the notion of
$n$-successful and local times. As usual we write $\log_{ 2} n$ for $\log\,\log n$.
\begin{lemma}\label{plem1} Let 
\[\mathcal{S}_n=\{x\in U_{ n}\,|\,\mbox{$x$ is
$n$-successful}\}.\] Then for some $N( \om)<\ff$ a.s., for all $n\geq N( \om)$ and
all $x\in \mathcal{S}_n$
\[
\frac{L^{ x}_{ T_{D(0, K_n)^{ c}}}}{ (\log K_n)^2}\geq  2a/\pi-2/\log_{ 2}
n.\]
\end{lemma}

\noindent {\bf Proof of Lemma~\ref{plem1}:} Recall that if $x$ is n-successful then
$N^x_{n,n}\geq \frak{n}_{n}-n=a(3n^2\log n)-n$.  Let $L^{ x,j}$ denote the
number of visits to $x$ during the
$j$'th excursion from $\partial D(x,r_{n,n})$ to $\partial D(x, r_{n,n-1})$. Then
 for any $0<\la<\ff$
\begin{eqnarray} P_x &:=&\PPP\(L^{ x}_{ T_{D(0, K_n)^{ c}}}
\leq (2a/\pi-2/\log_{ 2} n)(\log K_n)^2
\,, \,x\in \mathcal{S}_n\)
\nn\\ &\leq &\PPP\(\sum_{j=0}^{\frak{n}_{n}-n}L^{ x,j}
\leq (2a/\pi-1/\log_{ 2} n)(3n\log n)^2 
\)\nn\\ & \leq & \exp\(\la (2a/\pi-1/\log_{ 2} n)(3n\log n)^2\)
 E\(e^{ -\la \sum_{j=0}^{\frak{n}_{n}-n}L^{ x,j}} \).\label{pc3.6}
\end{eqnarray} 
 If $\tau$ denotes the first time that the $\frak{n}_{n}-n$'th excursion from 
$\partial D(x, r_{n,n-1})$ reaches
$\partial D(x,r_{n,n})$ then by the strong Markov property
\begin{eqnarray} &&\qquad E\(e^{ -\la \sum_{j=0}^{\frak{n}_{n}-n}L^{
x,j}} \)\label{pc3.6a}\\ &&=E\(e^{ -\la
\sum_{j=0}^{\frak{n}_{n}-n-1}L^{ x,j}}E^{X_{\tau }}\(e^{
-\la L^{ x}_{ T_{D(x, r_{n,n-1})^{ c}}}}\) \).\nn
\end{eqnarray} 
Set $\la=\phi/G_{D(x, r_{n,n-1}) }(x ,x)$.  
By (\ref{babe.1}), with $r=r_{n,n}=e^{ n}, R=r_{n,n-1}=n^{3 }e^{ n}$, 
for any 
$0<\phi\leq 1$ and large $n$
\begin{equation}
\sup_{y\in  \partial D(x,r_{n,n})}E^{y}\(e^{ -\la L^{ x}_{ T_{D(x, r_{n,n-1})^{
c}}}}\)\leq \exp\(-{\(1-1/2\log n\)\phi \over 1+\phi}\,3( \log n)/n
\).\label{pc3.6b}
\end{equation} 
Hence by induction
\begin{equation} E\(e^{ -\la \sum_{j=0}^{\frak{n}_{n}-n} L^{ x,j} } \)
\leq \exp\(-{( 1-1/\log n)\phi \over 1+\phi}9an( \log n)^{ 2}
\).\label{pc3.6c}
\end{equation} 
Then with this choice of $\la$, noting that $G_{D(x, r_{n,n-1})
}(x ,x)\sim { 2\over
\pi}n$ by (\ref{p2.5}), we have 
\begin{equation} P_x\leq \inf_{\phi>0 }\exp\(\Big\{\phi( 1-1/2\log_{ 2} n) - {(
1-1/\log n)\phi \over 1+\phi}\Big\}9an( \log n)^{ 2}\).\label{pc3.6d}
\end{equation} A straightforward computation shows that
\begin{equation}
\inf_{\varphi>0 }\(\varphi\al -{ \varphi\over 1+\varphi
}\bb\)=-\(\sqrt{\bb}-\sqrt{\al}\)^{ 2}\label{8.21}
\end{equation} which is achieved for $\varphi =\sqrt{\bb}/\sqrt{\al}-1$. Using this
in (\ref{pc3.6d}) we find that
\begin{equation} P_x\leq \exp\(-cn(\log n/\log_{ 2} n)^{ 2}\).\label{pc3.6e}
\end{equation} Note that $|U_{ n}|\leq e^{ cn\log n}$.
 Summing over all $x\in U_{ n}$ and then over $n$ and applying Borel-Cantelli will
then complete the proof of Lemma~\ref{plem1}.
\qed

Using this Lemma we see that to prove (\ref{q.2}) it will suffice to show that for any
$\de>0$ we can find $p_{ 0}>0$ and $N_{ 0}<\ff$ such that
\begin{equation}
\PPP\(\sum_{x\in U_{ n}}Y(n,x)\geq K_{ n}^{2-a-\de}\)
\geq p_{0}\label{pm3.6}
\end{equation} for all $n\geq N_{ 0}$. 

We will see by (\ref{pmomlb.1}) of the next Lemma below that for some sequence
 $\de_n\rar 0$ and constant $c>0$
\begin{equation}
\E\(\sum_{x\in U_{ n}}Y(n,x)\)=\sum_{x\in U_{ n}}
\bar{q}_{n,x}\geq cK_n^{2-a-\de_n}.\label{pm3.7}
\end{equation} Recall the Paley-Zygmund inequality (see
\cite[page 8]{Kahane}): for any $W\in L^{ 2}(
\Om)$ and $0<\la<1$
\begin{equation}
 \PPP(W\geq \la \E( W) )\geq ( 1-\la)^{ 2}{ (\E( W))^{ 2}\over
\E( W^{ 2})}.\label{PZ}
\end{equation} We will apply this with $W=W_{  n}=\sum_{x\in U_{  n}}Y(n,x)$.
Thus to complete the proof of (\ref{pm3.6}) it suffices to show that
\begin{equation}
\E\(\lc\sum_{x\in U_{ n}}Y(n,x)\rc^2\)\leq c\lc
\E\(\sum_{x\in U_{ n}}Y(n,x)\)\rc^2\label{p.2}
\end{equation} for some $c<\ff$ all $n$ sufficiently large. Furthermore, using
(\ref{pm3.7}) it suffices to show that
\begin{equation}
\E\(\sum_{\stackrel{x,y\in U_{ n}}{x\neq y}}Y(n,x)Y(n,y)\)\leq c\lc
\E\(\sum_{x\in U_{ n}}Y(n,x)\)\rc^2.\label{p.3}
\end{equation}

The next lemma, which provides estimates for the first and second moments of
$Y(n,x)$, will be proven in the following sections.

\begin{lemma}\label{pmomlb}  There exists $\de_n \to 0$ such that for all $n \geq
1$,
\begin{equation}\qquad
\bar{q}_{n,x}\geq Q_{ n}=\inf_{ x\in U_{ n}} \PPP(x\,
\mbox{ is $n$-successful}) \geq K_{n}^{-(a+\de_n)},
\label{pmomlb.1}
\end{equation}  and
\begin{equation} Q_n\geq c \bar{q}_{n,x}
\label{pmomlb.1j}
\end{equation} for some $c>0$ and all
$n$ and $x\in U_{ n}$.

There exists $C<\ff$ and $\de'_n \to 0$ such that for all $n$, $x \neq y$ and
$l(x,y)=\min \{m\,:
\,D(x,r_{n,m})\cap D(y,r_{n,m})=\emptyset\} \leq n$
\begin{equation}
\E(Y(n,x)Y(n,y))\leq C Q_n^2  ( l(x,y)!)^{3a+\de'_{l(x,y) }}\;.
\label{pm3.5}
\end{equation}  
\end{lemma}

In the sequel, we let $C_m$ denote generic finite constants that are independent of
$n$. The definition of $l(x,y)
\geq 1$ implies that
$|x-y|\leq 2r_{n,l(x,y)-1}$. Recall that there are at most
$C_0 r_{ n,l-1}^2\leq C_0K^{ 2}_n\,l^{ 6}( l!)^{ -6}$
 points 
$y$ in the ball of radius $2r_{n, l-1}$ centered at $x$. Thus, it follows from Lemma
\ref{pmomlb} that
\begin{eqnarray} &&
\sum_{\stackrel{x,y\in U_{ n}}{2r_{ n,n}\leq |x-y|\leq 2r_{
n,0}}}\E\(Y(n,x)Y(n,y)\)\label{p.5}\\ &&
  \leq C_1 \sum_{\stackrel{x,y\in U_{ n}}{2r_{ n,n}\leq |x-y|\leq 2r_{ n,0}}} Q_n^2 
( l(x,y)!)^{3a+\de'_{l(x,y) }}\nn\\ &&
\leq C_2 Q_n^2 K^{ 2}_n \sum_{l=1}^{n}  K^{ 2}_n\,l^{ 6}( l!)^{ -6}(
l!)^{3a+\de'_{l }}\nn\\  &&
\leq C_3 (K^{ 2}_n Q_n)^2  \sum_{l=1}^{n}  l^{ 6}(l!)^{-3(2-a)+\de'_{l }}\nn\\ 
 &&\leq C_4 (K^{ 2}_n Q_n)^2\leq C_5\lc \E\(\sum_{x\in U_{ n}}Y(n,x)\)\rc^2
\nn
\end{eqnarray} where we used the fact from (\ref{pmomlb.1}) that
\begin{equation}
K^{ 2}_n Q_n\leq c \sum_{x\in U_{ n}}\bar{q}_{n,x} =c\E\(\sum_{x\in U_{
n}}Y(n,x)\).\label{p.5f}
\end{equation}
Similarly
\begin{eqnarray} &&
\sum_{\stackrel{x,y\in U_{ n}}{|x-y|\leq 2r_{ n,n}}}\E\(Y(n,x)Y(n,y)\)
\leq \sum_{\stackrel{x,y\in U_{ n}}{|x-y|\leq 2r_{n,
n}}}\E\(Y(n,y)\)\label{p.5}\\&&
  \leq C_6 \sum_{x\,;\,|x|\leq 2r_{n, n}} K^{ 2}_nQ_n
  \leq C_7 \,e^{ 2n} K^{ 2}_nQ_n\leq C_8\lc \E\(\sum_{x\in
U_{ n}}Y(n,x)\)\rc^2
\nn
\end{eqnarray}
using (\ref{p.5f}) and (\ref{pm3.7}).
 This completes the proof of (\ref{p.3}) and hence of (\ref{p.1}).
\qed

\section{First moment estimates}\label{51momest} 

By (\ref{8.15fc}), (\ref{p2.5}) and (\ref{p1.2}) we have that for any $x\in
D(0,n)$
\begin{eqnarray}
 \PPP^x\(T_{0}<T_{D(0,n)^c}\) &=& {(2/ \pi)\log (n/ |x|)+O( |x|^{-1})
\over (2/ \pi)\log n +\ga +O(n^{ -1})} \nn\\
  &=& {\log (n/ |x|)+O( |x|^{-1}) \over
\log (n)}( 1+O( (\log n)^{-1}).  \label{p1.2c} 
\end{eqnarray}

By Exercise 1.6.8 of \cite{L1} we have that uniformly in $r<|x|<R$
\begin{equation}
\PPP^x\(T_{D(0,R)^c}<T_{ D(0,r)}\)= {\log (|x|/r)+O( r^{-1})\over
\log (R/ r)}\label{p1.2bfa}
\end{equation} and
\begin{equation}
\PPP^x\(T_{ D(0,r)}<T_{D(0,R)^c}\)= {\log (R/ |x|)+O( r^{-1})\over
\log (R/ r)}.\label{p1.2b}
\end{equation} 

\noindent {\bf Proof of Lemma \ref{pmomlb}}: For $x\in U_{ n}$ we begin by
getting bounds on the probability of reaching $\partial D(x,r_{ n,0})$ before 
$T_{\partial D(0,K_{ n})}$. Since 
\begin{eqnarray}
 \PPP \(T_{ \partial D(x,r_{n, 0})}<T_{\partial D(0,K_{ n})}\) 
\geq  \PPP \(T_{ \partial D(x,r_{ n,0})}<T_{\partial D(x,{1 \over 2}K_{ n})}
\)\label{pr.1}  
\end{eqnarray}  we see from (\ref{p1.2b}) that uniformly in $n$ and $x\in U_{ n}$
\begin{equation}
\PPP \(T_{\partial  D(x,r_{n, 0})}<T_{\partial D(0,K_{ n})}\) \geq c\label{pr.2}
\end{equation}  for some $c>0$. And since for $x\in U_{ n}$ and $y\in
\partial D(x,r_{n,0})$
\begin{eqnarray}
&&
\PPP^{ y} \(T_{\partial  D(x,r_{n, 1})}<T_{\partial D(x,{1 \over 2}K_{
n})}\)\label{pr.3}\\ &&\hspace{ .3in}\leq
 \PPP^{ y} \(T_{ \partial D(x,r_{n, 1})}<T_{\partial D(0,K_{ n})}\)\nn\\
&&\hspace{ .6in}\leq  \PPP^{ y} \(T_{\partial  D(x,r_{n, 1})}<T_{\partial D(x,2K_{
n})}\)
\nn
\end{eqnarray}  we see from (\ref{p1.2b}) that uniformly in $n$, $x\in U_{ n}$
and $y\in
\partial D(x,r_{n,0})$
\begin{equation}
 c/\log n\leq \PPP^{ y} \(T_{ \partial D(x,r_{n, 1})}<T_{\partial D(0,K_{ n})}\) \leq
c'/\log n.\label{pr.4}
\end{equation} 
Similarly, since for $x\in U_{ n}$ and $y\in
\partial D(x,r_{n,0})$
\begin{eqnarray}
 \PPP^{ y} \(T_{\partial D(0,K_{ n})}<T_{ \partial D(x,r_{n, 1})}\) 
\geq   \PPP^{ y} \(T_{\partial D(x,2K_{ n})}<T_{ \partial D(x,r_{n, 1})}\)
\label{pr.3}
\end{eqnarray}
we see from (\ref{p1.2bfa}) that uniformly in $n$, $x\in U_{ n}$
and $y\in
\partial D(x,r_{n,0})$
\begin{equation}
 \PPP^{ y} \(T_{\partial D(0,K_{ n})}<T_{ \partial D(x,r_{n, 1})}\)  \geq
c>0. 
\label{pr.5}
\end{equation}

These bounds will be used for excursions at the `top' levels. To bound excursions at 
`intermediate' levels we note that  using (\ref{p1.2bfa}), we have uniformly for
$x\in
\partial D(0,r_{n,l})$, with 
$1\leq l\leq n-1$
\begin{eqnarray}
 &&
 \PPP^x\(T_{\partial D(0,r_{n, l-1})}<T_{ \partial D(0,r_{ n,l+1})}\) 
 =1/2+O(n^{ -8})\label{ps.3}
\end{eqnarray}  
and consequently also
\begin{eqnarray}
 \PPP^x\(T_{ \partial D(0,r_{ n,l+1})}<T_{\partial D(0,r_{n, l-1})}\) \label{ps.5}
 =1/2+O( n^{ -8}). \label{ps.5}
\end{eqnarray}

Let $\bar{m}=(m_{ 2}, m_{ 3},\ldots, m_{ n})$
 and set $|\bar{m}|=2\sum_{ j= 2}^{ n}m_{ j}+1$. Let $\mathcal{H}_{ n}(\bar{m})$,
 be the collection of maps, (`histories'),
\[\varphi: \{0,1,\ldots,  |\bar{m}|  \}\mapsto \{0,1,\ldots,n \}\] such that
$\varphi(0)=1,\,\varphi( j+1)=\varphi( j)\pm 1,\,|\bar{m}|
=\inf\{j\,;\,\varphi(j)=0
\}$
 and the number of upcrossings from $\ell-1$ to
$\ell$
\[u( \ell)=:|\{( j,j+1) \,|\, ( \varphi(j),\varphi( j+1))=(\ell-1,\ell)  \}|.\]

The number of ways to partition the $u( \ell+1)$  upcrossings from $\ell$ to
$\ell+1$ among and after the $u( \ell)$  upcrossings from $\ell-1$ to
$\ell$ is
\begin{equation} { u( \ell+1)+u( \ell)-1\choose u( \ell)-1}.\label{44.1}
\end{equation}
 Since $u( \ell)=m_\ell$ and the mapping $\varphi$ is completely determined once
 we know the relative order of all its upcrossings
\begin{equation} |\mathcal{H}_{ n}(\bar{m})|=
\prod_{\ell= 2}^{ n-1}{m_{ \ell+1}+m_{ \ell}-1\choose m_{
\ell}-1}.\label{44.2} 
\end{equation}

 To each random walk path we
assign a `history' 
$h(\om )$ as follows.
 Let $\tau( 0)$ be the time of the first visit to 
$\partial D(x,r_{n, 1})$, and  define $\tau( 1), \tau( 2),
\ldots$ to be the successive hitting times of different elements of
\[\{\partial D(x,r_{n,0}),
\ldots,
\partial D(x,r_{n, n})\}\] until the first downcrossing from 
$\partial D(x,r_{ n,1})$ to
$\partial D(x,r_{ n,0})$. Setting $\Phi( y)=k$ if $y\in \partial D(x,r_{ n,k})$,
 let
$h(\om )( j)=\Phi( \om (\tau( j) ))$.  Let $h_{|_{k} }$ be the restriction of $h$ to 
$\{0,\ldots,k  \}$.  We claim that uniformly for any
 $\varphi\in \mathcal{H}_{ n}(\bar{m})$ and $z\in \partial D(x,r_{n,1})$
\begin{eqnarray}
 &&
  P^{ z}\lc h_{|_{ |\bar{m}|} }=\varphi\rc 
 =\( { 1\over 2} \)^{|\bar{m}|-m_{ n} }\lc 1+O(n^{ -8})
 \rc^{|\bar{m}| -m_{ n} }.\label{44.3}
\end{eqnarray} 
To see this, simply use the Markov property successively at the
times 
\[\tau( 0), \tau( 1),
\ldots, \tau(|\bar{m}|-1)\] and then use (\ref{ps.3}), (\ref{ps.5}). (The $m_{ n}$
 downcrossings from $n$ to $n-1$ come `for free').

Writing $m\stackrel{k}{\sim} \frak n_{k}$ if $m=1$ for $k<k_{ 0}$ and 
$|m-\frak n_k| \leq k$ for $k\geq k_{ 0}$ we see 
 that uniformly in $m_{ n}\stackrel{ n}{\sim}\frak n_{ n}$ we have that
$\lc 1+O(n^{ -8})
 \rc^{|\bar{m}| }=1+O(n^{-4})$. Combining this with (\ref{44.2}) and (\ref{44.3})
we see that  uniformly in $z\in \partial D(x,r_{n,1})$
\begin{eqnarray} &&
\sum_{\stackrel{m_{ 2}, \ldots, m_{ n}} {m_{\ell}\stackrel{\ell}{\sim}\frak n_{
\ell}} }
 P^{ z}\lc h_{|_{ |\bar{m}|} }\in
\mathcal{H}_{ n} (\bar{m})\rc\label{44.5m}\\ &&
 =( 1+O(n^{-4}))\sum_{\stackrel{m_{ 2}, \ldots, m_{ n}}
{m_{\ell}\stackrel{\ell}{\sim}\frak n_{
\ell}} }\( { 1\over 2} \)^{|\bar{m}|-m_{ n} }
\prod_{\ell= 2}^{ n-1}{m_{ \ell+1}+m_{ \ell}-1\choose m_{
\ell}-1} \nonumber\\ &&
 =( 1+O(n^{-4})){ 1\over 4}\sum_{\stackrel{m_{ 2}, \ldots, m_{ n}}
{m_{\ell}\stackrel{\ell}{\sim}\frak n_{
\ell}} }
\prod_{\ell= 2}^{ n-1}{m_{ \ell+1}+m_{ \ell}-1\choose m_{
\ell}-1}\( { 1\over 2} \)^{m_{ \ell+1}+m_{ \ell} }. \nonumber
\end{eqnarray}
Here we used the fact that $m_{ 2}=1$ so that $|\bar{m}|-m_{ n} =2+\sum_{\ell=
2}^{ n-1}(m_{ \ell+1}+m_{ \ell})$.

\begin{lemma}\label{lem-stirling}
For some $C=C(a) <\infty$ and
all $k\geq 2$,
$|m-\frak n_{k+1}| \leq k+1$, $|\ell+1-\frak n_{k}| \leq  k$,
\begin{equation} { C^{ -1}k^{-3a-1}\over \sqrt{\log k}}
\leq {m+\ell\choose \ell}\( { 1\over 2} \)^{m+\ell+1}
\leq { Ck^{-3a-1}\over \sqrt{\log k}}.\label{44.5}
\end{equation}
\end{lemma}

\noindent
{\bf Proof of Lemma \ref{lem-stirling}}:
It suffices to consider $k \gg 1$
in which case the binomial coefficient in (\ref{44.5})
is well approximated by Stirling's formula
\[
m!=\sqrt{2\pi}m^m e^{-m}\sqrt{m} (1+o(1)) \;.
\]
With $\frak n_k = 3a k^2 \log k$ it follows that for some $C_1<\ff$
and all $k$ large enough, if $|m-\frak n_{k+1}| \leq 2 k$,
$|\ell-\frak n_{k}| \leq 2 k$ then
\begin{equation}
|\frac{m}{\ell} - 1 - \frac{2}{k}| \leq \frac{C_1}{k \log k} \;.
\label{m1.3la}
\end{equation}
Hereafter, we use the notation
 $f \sim g$ if $f/g$ is bounded and
bounded away from zero as $k \to \infty$, uniformly in
$\{m:\,|m-\frak n_{k+1}| \leq 2 k\}$ and $\{\ell:|\ell-\frak n_{k}| \leq 2 k\}$.
We then have by the preceeding observations that
\begin{equation}
{m+\ell\choose \ell}\( { 1\over 2} \)^{m+\ell+1} \sim
\frac{(m+\ell)^{m+\ell}}{\sqrt{\ell}\, \ell^\ell m^m}
\( { 1\over 2} \)^{m+\ell} \sim \frac{\exp(-\ell I(
{m \over \ell}))}{\sqrt{k^2 \log k}} \;,
\label{Lkest}
\end{equation}
where
$$
I(\la)= - (1+\la) \log (1+\la) + \la \log \la + \la \log 2 + \log 2 \;.
$$
The function $I(\la)$ and its
first order derivative vanishes at $1$, with
the second derivative $I_{\la \la} (1) = 1/2$.
Thus, by a Taylor expansion to second order of $I(\la)$ at $1$,
the estimate (\ref{m1.3la}) results with
\begin{equation}
| I({m \over \ell}) - \frac{1}{k^2} | \leq \frac{C_2}{k^2 \log k}
\label{Ipest}
\end{equation}
for some $C_2 < \infty$, all $k$ large enough and $m,\ell$ in the range
considered here.  Since
$|\ell - 3a k^2 \log k| \leq 2k$,
combining (\ref{Lkest}) and (\ref{Ipest}) we establish (\ref{44.5}).
\qed

Using the last Lemma we have that  
\begin{eqnarray} &&
\sum_{\stackrel{m_{ 2}, \ldots, m_{ n}} {m_{\ell}\stackrel{\ell}{\sim}\frak n_{
\ell}} }
\prod_{\ell= 2}^{ n-1}{ C^{ -1}\ell^{-3a-1}\over \sqrt{\log \ell}}\nn\\ &&
\hspace{ .5in}\leq
\sum_{\stackrel{m_{ 2}, \ldots, m_{ n}} {m_{\ell}\stackrel{\ell}{\sim}\frak n_{
\ell}} }
\prod_{\ell= 2}^{ n-1}{m_{ \ell+1}+m_{ \ell}-1\choose m_{
\ell}-1}\( { 1\over 2} \)^{m_{ \ell+1}+m_{ \ell} }\label{44.5n}\\ && 
\hspace{ 1in}\leq   
\sum_{\stackrel{m_{ 2}, \ldots, m_{ n}} {m_{\ell}\stackrel{\ell}{\sim}\frak n_{
\ell}} }
\prod_{\ell= 2}^{ n-1}{ C\ell^{-3a-1}\over \sqrt{\log \ell}}.  \nonumber
\end{eqnarray}
Using the fact that $|\{ m_{\ell}\,|\,m_{\ell}\stackrel{\ell}{\sim}\frak n_{
\ell}\}|=2\ell+1$, this shows that for some $C_1<\ff$,
\begin{eqnarray} &&
n\prod_{\ell= 2}^{ n-1}{ C_1^{ -1}\ell^{-3a}\over \sqrt{\log \ell}}\nn\\ &&
\hspace{ .5in}\leq
\sum_{\stackrel{m_{ 2}, \ldots, m_{ n}} {m_{\ell}\stackrel{\ell}{\sim}\frak n_{
\ell}} }
\prod_{\ell= 2}^{ n-1}{m_{ \ell+1}+m_{ \ell}-1\choose m_{
\ell}-1}\( { 1\over 2} \)^{m_{ \ell+1}+m_{ \ell} }\label{44.5p}\\ && 
\hspace{ 1in}\leq   
n\prod_{\ell= 2}^{ n-1}{ C_1\ell^{-3a}\over \sqrt{\log \ell}}.  \nonumber
\end{eqnarray}

Since for any $c<\ff$, for some $\zeta_{ n},\zeta'_{ n}\rar 0$ 
\begin{equation}
nc^{ n}\prod_{\ell= 2}^{ n-1}\log \ell=n^{n\zeta_{ n} }=( n!)^{\zeta'_{ n} }
\label{44.5q}
\end{equation}
we see that  for some $\de_{1, n}, \de_{ 2,n}\rar 0$
\begin{equation}
\sum_{\stackrel{m_{ 2}, \ldots, m_{ n}} {m_{\ell}\stackrel{\ell}{\sim}\frak n_{
\ell}} }
\prod_{\ell= 2}^{ n-1}{m_{ \ell+1}+m_{ \ell}-1\choose m_{
\ell}-1}\( { 1\over 2} \)^{m_{ \ell+1}+m_{ \ell} }
 =( n!)^{-3a-\de_{1, n}}=r_{n,0}^{-a-\de_{2, n}}.\label{44.6}
\end{equation}

(\ref{pr.2})-(\ref{pr.5}) and (\ref{44.5m}) show that for some $0<c,c'<\ff$
\begin{eqnarray} &&
{c \over \log n}\sum_{\stackrel{m_{ 2}, \ldots, m_{ n}}
{m_{\ell}\stackrel{\ell}{\sim}\frak n_{
\ell}} }
\prod_{\ell= 2}^{ n-1}{m_{ \ell+1}+m_{ \ell}-1\choose m_{
\ell}-1}\( { 1\over 2} \)^{m_{ \ell+1}+m_{ \ell} }\label{44.6f}\\ && 
\hspace{ .3in}\leq Q_{ n}=\inf_{ x\in U_{ n}} \PPP(x\,
\mbox{ is $n$-successful})\nonumber\\ &&
\hspace{ .6in}\leq {c' \over \log n}
\sum_{\stackrel{m_{ 2}, \ldots, m_{ n}} {m_{\ell}\stackrel{\ell}{\sim}\frak n_{
\ell}} }
\prod_{\ell= 2}^{ n-1}{m_{ \ell+1}+m_{ \ell}-1\choose m_{
\ell}-1}\( { 1\over 2} \)^{m_{ \ell+1}+m_{ \ell} }.\nonumber
\end{eqnarray}

Together with (\ref{44.6}) this gives 
(\ref{pmomlb.1}) and (\ref{pmomlb.1j}).

 Let $N_{n,i,m,k}^x$ denote the number of excursions 
from
$ \partial D(x,r_{n,k-1})$ to $
\partial D(x,r_{n,k})$  until completion of the first
$m$ excursions from
$ \partial D(x,r_{n,i})$ to $
\partial D(x,r_{n,i-1})$. We note here for later reference that the analysis of this
section shows that uniformly in 
$m_{k}\stackrel{k}{\sim}\frak n_{k},\; k=l,l+1,\ldots,n$ and 
$z\in \partial D(x,r_{n,l})$
\begin{eqnarray}
 &&
 \PPP^{ z} \( N_{n,l,m_{ l},k}^x =m_{k},\; k=l+1,\ldots,n\) 
\label{44.8}\\
 &&
 = (1+O(n^{-4}))\prod_{ k=l}^{ n-1}{m_{ k+1}+m_{ k}-1\choose m_{ k}-1}\( { 1\over
2} \)^{m_{ k+1}+m_{ k} }.\nn
\end{eqnarray} 
Our analysis also shows that for some $\de_{3, l}\rar 0$
\begin{equation}
 \inf_{\stackrel{ m_{ l}} {m_{l}\stackrel{l}{\sim}\frak n_{ l}} 
}\left(\sum_{\stackrel{m_{ 2}, \ldots, m_{ l-1}} {m_{k}\stackrel{k}{\sim}\frak
n_{ k}} }\prod_{ k=2}^{ l-1}{m_{ k+1}+m_{ k}-1\choose m_{ k}-1}\( { 1\over 2}
\)^{m_{ k+1}+m_{ k} }\right)\geq ( ( l-1)!)^{ -3a-\de_{3, l}},\label{44.8y}
\end{equation}
and since 
\begin{eqnarray} &&
\sum_{\stackrel{m_{ 2}, \ldots, m_{ n}} {m_{\ell}\stackrel{\ell}{\sim}\frak n_{
\ell}} }
\prod_{\ell= 2}^{ n-1}{m_{ \ell+1}+m_{ \ell}-1\choose m_{
\ell}-1}\( { 1\over 2} \)^{m_{ \ell+1}+m_{ \ell} }\label{44.8z}\\
&&\geq\sum_{\stackrel{m_{ l}, \ldots, m_{ n}} {m_{k}\stackrel{k}{\sim}\frak n_{
k}} }\prod_{ k=l}^{ n-1}{m_{ k+1}+m_{ k}-1\choose m_{ k}-1}\( { 1\over
2} \)^{m_{ k+1}+m_{ k} }\nn\\ && 
 \inf_{\stackrel{ m_{ l}} {m_{l}\stackrel{l}{\sim}\frak n_{ l}} 
}\left(\sum_{\stackrel{m_{ 2}, \ldots, m_{ l-1}} {m_{k}\stackrel{k}{\sim}\frak
n_{ k}} }\prod_{ k=2}^{ l-1}{m_{ k+1}+m_{ k}-1\choose m_{ k}-1}\( { 1\over 2}
\)^{m_{ k+1}+m_{ k} }\right)\nonumber
\end{eqnarray}
(here we used the bound that for $C( i,j),D( i,j)$ non-negative, 
$\sum_{ i,j,k}C( i,j)D( j,k)=\sum_{j}\sum_{ i}C( i,j)\sum_{k}D(  j,k)
\geq \(\sum_{ i,j}C( i,j)\) \inf_{ j}\sum_{k}D(  j,k)$\,\,),
we see from (\ref{44.8}) and (\ref{44.6f}) that uniformly in 
$z\in \partial D(x,r_{n,l})$
\begin{equation}
 \sum_{\stackrel{m_{ l}, \ldots, m_{ n}} {m_{k}\stackrel{k}{\sim}\frak n_{ k}}
}\PPP^{ z}
\( N_{n,l,m_{ l},k}^x =m_{k},\; k=l+1,\ldots,n\) 
\leq C\log n\, Q_n  \,\,( l!)^{3a+\de_{3,l }}.\label{44.9s}
\end{equation}
Thus from (\ref{pr.4}) we have 
\begin{equation}
 \sum_{\stackrel{m_{ l}, \ldots, m_{ n}} {m_{k}\stackrel{k}{\sim}\frak n_{ k}}
}\PPP
\( N_{n,l,m_{ l},k}^x =m_{k},\; k=l+1,\ldots,n\) 
\leq C\, Q_n  \,\,( l!)^{3a+\de_{3,l }}.\label{44.9q}
\end{equation}

Similarly, uniformly in 
$m_{k}\stackrel{k}{\sim}\frak n_{k},\; k=2,3,\ldots,l$ and $z\in \partial
D(x,r_{n,1})$
\begin{eqnarray}
 &&
 \PPP^{ z} \( N_{n,1,1,k}^x =m_{k},\; k=2,\ldots,l\) 
\label{44.8m}\\
 &&
 = (1+O(n^{-4}))\prod_{ k=2}^{ l-1}{m_{ k+1}+m_{ k}-1\choose m_{ k}-1}\( { 1\over
2} \)^{m_{ k+1}+m_{ k} }.\nn
\end{eqnarray}

\section{Second moment estimates}\label{52momest}

We begin by defining the
$\si$-algebra $\GG_{n, l}^x$ of excursions from
$\partial D(x,r_{n, l-1})$ to $\partial D(x,r_{n, l})$. To this end, fix  
$x \in \Z^2$, let $\taup_0=0$ and for $i=1,2,\ldots$ define
\begin{eqnarray*}
\tau_{i} & = & \inf \{ t \geq \taup_{i-1}  :\; X_t \in \partial D(x,r_{n, l}) \} \,, \\
\taup_{i} & = & \inf \{ t \ge \tau_{i} :\; X_t \in  \partial D(x,r_{n, l-1})\}.
\end{eqnarray*} Then $\GG_{n, l}^x$ is the $\sigma$-algebra generated by the
excursions
$\{ e^{(j)}, j =1,\ldots \}$, where
$e^{(j)}=\{ X_{t} : \taup_{j-1}\leq t \leq \tau_{j} \}$ is the $j$-th excursion from
$\partial D(x,r_{ n,l-1})$ to $\partial D(x,r_{n, l})$ (so for $j=1$ we do begin at
$t=0$).

The following Lemma is proven in the next section. Note that for any 
$\si$-algebra $\GG$ and event $B\in\GG$, we have 
 $P( A,B\,|\, \GG)=P( A\,|\, \GG)1_{ \{ B\}}$.

\begin{lemma}[Decoupling Lemma]\label{cond-cor} Let $\Gamma_{
n,l}^y =\{ N^y_{n,i} = m_i ; i=l+1,l+2,\ldots,n \}$. Then, uniformly over all $l\leq n
,\,$ 
$m_l \stackrel{l}{\sim} \frak n_l$, $\{m_i: i=l+1,l+2,\ldots,n\}$,
$y \in U_{ n}$ and $x_0,x_1 \in \Z^2 \setminus D(y,r_{n,l})$,
\begin{eqnarray} &&
\PPP^{x_0} ( \Gamma_{ n,l}^y\, ,\,N_{n,l}^y=m_l\, |\,
{\mathcal G}_{ n,l}^y)\label{new1.3e}\\ &&
 = (1+O(n^{-1} \log n))  \PPP^{x_1} (\Gamma_{ n,l}^y 
\, |\,N_{n,l}^y=m_l)1_{ \{ N_{n,l}^y=m_l\}} \nonumber
\end{eqnarray}
\end{lemma}

\begin{remark}{\rm The intuition behind the Decoupling Lemma is that what
happens `deep inside' $D(y,r_{n, l})$, e.g. $\Gamma_{ n,l}^y $, is
`almost' independent of what happens outside $D(y,r_{n, l})$, i.e. ${\mathcal G}_{
n,l}^y$.}
\end{remark}

\noindent{\bf Proof of (\ref{pm3.5}):} Recall that $\frak{n}_k=3a k^2 \log k$ and
that we write 
$m\stackrel{k}{\sim} \frak{n}_{k}$ if $m=1$ for $k<k_{ 0}$ and 
$|m-\frak n_k| \leq k$ for $k\geq k_{ 0}$. Relying upon the first moment estimates
and Lemma \ref{cond-cor}, we next prove the second moment estimates
\req{pm3.5}. Take $x,y \in U_{n}$ with $l( x,y)=l-1$. Thus $2 r_{n,l-1}+2 \leq 
|x-y|< 2 r_{n,l-2}+2$ for some $2 \leq l \leq n$. Since
$r_{n,l-3} -  r_{n,l-2} \gg  2 r_{n,l-1}$, it is easy to see that
$D(y,r_{n,l-1}) \cap \partial D(x,r_{n,k}) =\emptyset$ for all $k \neq l-2$.
 Replacing hereafter $l$ by $l
\wedge (n-3)$, it follows that for
$k\neq l-1,l-2$,
 the events $\{ N_{n,k}^x \stackrel{k}{\sim} \frak{n}_{k} \}$ are measurable on the
$\sigma$-algebra $\GG^y_{n,l}$. With $J_l:=\{ l+1,\ldots,n \}$ and
$I_l:=\{ 2,\ldots,l-3,l,\ldots,n \}$, set 
$\wt{\Ga}^{ y}_{ n}(J_l )=\{ N_{n,k}^y\stackrel{k}{\sim} \frak{n}_{k},\, k \in
J_{l}\}$
 and 
$\wt{\Ga}^{ x}_{ n}(I_l )=\{ N_{n,k}^x \stackrel{k}{\sim} \frak{n}_{k},\; k \in
I_l\}$.  We note that
\begin{eqnarray}  &&
\{ x, y \, \mbox{ are $n$-successful} \}\label{jstar}\\ &&\hspace{ 1in}
\subseteq \bigcup_{m_l\stackrel{l}{\sim} \frak{n}_{l} }
 \wt{\Ga}^{ y}_{ n}(J_{ l} )\cap
\{ N^y_{n,l} = m_{l}\}\cap \wt{\Ga}^{ x}_{ n}(I_l ).\nn
\end{eqnarray}
  Applying \req{new1.3e}, we have that for some universal constant
$C_3<\infty$,
\begin{eqnarray} &&
\PPP\(x\, \mbox{ and } \, y \, \mbox { are $n$-successful} \)\label{52mom}\\ &&
\leq    \sum_{m_l\stackrel{l}{\sim} \frak{n}_{l}}
P \left[  \wt{\Ga}^{ y}_{ n}(J_{ l} ),
\{ N^y_{n,l} = m_{l}\}, \wt{\Ga}^{ x}_{ n}(I_l )
\right]\nonumber\\ &&
=   \sum_{m_l\stackrel{l}{\sim} \frak{n}_{l}}
\E \left[ \PPP ( \wt{\Ga}_{ n}^{ y}(J_{ l} ),\, N^y_{n,l} = m_{l} \,\big| \,
 {\mathcal G}^y_{n,l})  \,;\wt{\Ga}_{ n}^{ x}(I_l )
\right]\nonumber\\ &&\leq C_3
\PPP (\wt{\Ga}_{ n}^{ x}(I_l ))
\sum_{m_l\stackrel{l}{\sim} \frak{n}_{l}}
  \PPP ( \wt{\Ga}_{ n}^{ y}(J_{ l} )\,
\big| \, N^y_{n,l} = m_{l} )\nn
\end{eqnarray}
Using the Markov property and (\ref{44.9q}), for some universal constant
$C_5<\infty$,
\begin{eqnarray}\qquad
\sum_{ m_{l}\stackrel{l}{\sim} \frak{n}_{l}}\PPP (
 \wt{\Ga}_{ n}^{ y}(J_{ l} ) \, \big|
   \, N^y_{n,l} = m_{l} )
 \leq C_5  Q_n\,\,  ( l!)^{3a+\de_{3,l }}.\label{244.9} 
\end{eqnarray}  Similarly, with 
$M_{l-3 } :=\{2,\ldots,l-3 \}$, and
\[\wt{\Ga}_{ n}^{ x}(M_{l-3 })=\{ N_{n,k}^x \stackrel{k}{\sim} \frak{n}_{k},\; k \in
M_{l-3 }
\},\]  
\req{new1.3e} shows that
\begin{eqnarray} \qquad
\PPP \( \wt{\Ga}_{ n}^{ x}(I_l )\)
&\leq &
\sum_{m_{l}\stackrel{l}{\sim} \frak{n}_{l}}
\E \left[ \PPP (\wt{\Ga}_{ n}^{ x}(J_{ l} ),\, N^x_{n,l} = m_{l} \,
\big|
  {\mathcal G}^x_{n,l} ) \,;
\wt{\Ga}^{ x}(M_{l-3 }) \right]\label{244.7}\\
   &\leq & C_6 \PPP \( \wt{\Ga}_{ n}^{ x}(M_{l-3 } )\)
\sum_{m_{l}\stackrel{l}{\sim} \frak{n}_{l}}
\PPP (\wt{\Ga}_{ n}^{ x}(J_{ l} ) \, \big|
\, N^x_{n,l} = m_{l}).\nn
\end{eqnarray}  
Using (\ref{44.5m}), (\ref{44.6}) and  (\ref{244.9}) we get that, for some 
$\de_{4,l}\rar 0$
\begin{equation}
\PPP \(\wt{\Ga}_{ n}^{ x}(I_l ) \)
\leq C_7 l^{ 15}\,  ( l!)^{\de_{4,l }}\,Q_n.\label{244.10}
\end{equation} Putting (\ref{52mom}), (\ref{244.9}) and (\ref{244.10}) together
and adjusting $C$ and $\de'_{ l-1}$ proves (\ref{pm3.5}) for $l( x,y)=l-1$. 
\qed

\section{Harnack inequalities and approximate decoupling}\label{harn}

The goal of this section is to prove the Decoupling Lemma, Lemma \ref{cond-cor}.
Since what
happens `deep inside' $D(y,r_{n, l})$, e.g. $\Gamma_{ n,l}^y $,
depends on what happens outside $D(y,r_{n, l})$, i.e. on ${\mathcal G}_{
n,l}^y$, only through the initial and end points of the excursions from 
$\partial D(x,r_{ n,l})$ to $\partial D(x,r_{n, l-1})$, we begin by studying the
dependence on these initial and end points. The following Harnack
inequality plays a crucial role in this analysis.

Define the hitting distribution of $A$ by
\begin{equation} H_{A}(x,y)=P^x(X_{T_{A}}=y).\label{p1.4}
\end{equation}

\begin{lemma}\label{lem-H} Uniformly for $x,x'\in D(0,\varepsilon n)$,
$\varepsilon <1/4$ and $y\in D(0,n)^c$
\begin{eqnarray}
 H_{D(0,n)^c}(x,y)=\(1+O(\varepsilon )\)
H_{D(0,n)^c}(x',y).\label{p1.5}
\end{eqnarray} 

Furthermore, if $\varepsilon '<\varepsilon $ are such that
\[
\inf_{ x\in \partial D(0,\varepsilon  n)}
\PPP^{ x}(T_{D(0,n)^c }<T_{ \partial D(0,\varepsilon' n) })\geq 1/4,
\] 
 then uniformly in $x\in \partial D(0,\varepsilon  n)$
and
$y\in D(0,n)^c$,
\bea &&
\PPP^{ x}(X_{ T_{D(0,n)^c }}=y\,,\,T_{D(0,n)^c }<T_{\partial  D(0,\varepsilon' n) } )
\label{p1.5h}\\ && =\(1+O(\varepsilon )\)
\PPP^{ x}(T_{D(0,n)^c }<T_{\partial  D(0,\varepsilon ' n) })H_{D(0,n)^c}(x,y).\nn
\eea
\end{lemma}

\noindent {\bf Proof of Lemma~\ref{lem-H}:} 
(\ref{p1.5}) is formula (2.7) of \cite{L1}.

Turning to \req{p1.5h}, we have
\bea &&
\label{y1.5h}
\PPP^{ x}(X_{ T_{D(0,n)^c }}=y\,,\,T_{D(0,n)^c }<T_{ \partial  D(0,\varepsilon ' n)
}  )
\\ && =H_{ D(0,n)^c}(x,y)-
\PPP^{ x}(X_{ T_{D(0,n)^c }}=y\,,\,T_{D(0,n)^c }>T_{ \partial  D(0,\varepsilon ' n)
}).\nn
\eea By the strong Markov property at $T_{\partial  D(0,\varepsilon ' n)}$
\bea &&
\PPP^{ x}(X_{ T_{D(0,n)^c }}=y\,,\,T_{D(0,n)^c }>T_{\partial   D(0,\varepsilon ' n) }
)\label{y1.7h}\\ && =\E^{ x}(H_{D(0,n)^c}(X_{T_{ \partial  D(0,\varepsilon ' n)
}},y) ;\,T_{D(0,n)^c  }>T_{\partial   D(0,\varepsilon ' n) }).\nn
\eea By (\ref{p1.5}), uniformly in $w\in \partial  D(0,\varepsilon ' n)$,
\[ H_{ D(0,n)^c}(w,y) =
\(1+O(\varepsilon )\) H_{ D(0,n)^c}(x,y))\,.
\] Substituting back into (\ref{y1.7h}) we have
\beaa &&
\PPP^{ x}(X_{ T_{D(0,n)^c  }}=y\,,\,T_{D(0,n)^c  }>T_{\partial   D(0,\varepsilon '
n) } )\label{1.5tight}
\\ && =\(1+O(\varepsilon )\)
\PPP^{ x}(T_{D(0,n)^c  }>T_{\partial   D(0,\varepsilon ' n)
})H_{D(0,n)^c}(x,y).\nn
\eeaa

Combining this with (\ref{y1.5h}) we obtain
\bea &&
\PPP^{ x}(X_{ T_{D(0,n)^c  }}=y\,,\,T_{D(0,n)^c  }<T_{\partial  D(0,\varepsilon ' n)
} )\label{1.5tighter}
\\ && =\(\PPP^{ x}(T_{D(0,n)^c  }<T_{ \partial  D(0,\varepsilon ' n)
})+O(\varepsilon )\) H_{D(0,n)^c}(x,y).\nn
\eea Using the assumptions of the lemma  we obtain \req{p1.5h} which completes
the proof of Lemma~\ref{lem-H}.                  
\qed

 Consider a random path beginning at 
$z \in  \partial D(0,r_{ n,l})$. We will show that for
$n$ large, a certain
$\si$-algebra of excursions of the path from $ \partial D(0,r_{n, l+1})$ to 
$\partial D(0,r_{n,  l})$ prior to $T_{\partial D(0,r_{n, l-1})}$, is almost
independent of the choice of initial point
$z \in  \partial D(0,r_{n, l})$ and  final point $w \in   \partial D(0,r_{n,
l-1})$.  Let $\tau_0=0$ and for $i=0,1,\ldots$ define
\begin{eqnarray*}
\tau_{2i+1} & = & \inf \{ t \geq \tau_{2i}  :\; X_t \in  \partial D(0,r_{ n, l+1}) \cup 
\partial D(0,r_{ n,l-1})\} \\
\tau_{2i+2} & = &
   \inf \{ t \ge \tau_{2i+1} :\; X_t \in  \partial D(0,r_{n, l})\}\,.
\end{eqnarray*} Abbreviating $\bar{\tau}=T_{ \partial D(0,r_{n, l-1})}$ note that
$\bar{\tau}=\tau_{2I+1}$ for some (unique) non-negative integer $I$.
As usual, $\FF_{j}$ will denote the $\si-$algebra generated by $\{X_{
l},\,l=0,1,\ldots, j \}$, and for any stopping time $\tau$,
$\FF_{\tau}$ will denote the collection of events $A$ such that 
$A\cap \{ \tau=j\}\in \FF_{j}$ for all $j$.

Let $\HH_{n, l}$ denote the $\si$-algebra generated by the excursions of the path
from $ \partial D(0,r_{n,  l+1})$ to  $\partial D(0,r_{ n,l})$, prior to $T_{\partial 
D(0,r_{n, l-1}) }$. Then
$\HH_{n, l}$ is the $\sigma$-algebra generated by the excursions
$\{ v^{(j)}, j =1,\ldots, I \}$, where
$v^{(j)}=\{ X_{t} : \tau_{2j-1}\leq t \leq \tau_{2j} \}$ is the $j$-th excursion  from
$ \partial D(0,r_{ n,l+1})$ to  $\partial D(0,r_{ n,l})$. We denote by 
$\mathcal{U}(\{ \HH_{n, l}\})$ the collection of sequences of sets
$\{B_{ n}\,; n=1,2,\ldots \},\,B_{ n}\in\HH_{ n,l}$ such that uniformly in $x,x'\in
\partial D(0,r_{n, l+1})$
\begin{equation}
\PPP^{ x} (B_{ n}) = (1 +O( n^{ -4})) \PPP^{ x'} (B_{ n})
\,.
\label{m1.5unif}
\end{equation}

\begin{lemma}\label{lemprobends} For a random walk path starting at 
$z \in  \partial D(0,r_{ n,l})$, let $\HH_{ n,l}$ denote the $\si$-algebra
generated by the excursions of the path from $ \partial D(0,r_{n, l+1})$ to 
$\partial D(0,r_{n, l})$,
 prior to
$T_{\partial  D(0,r_{n,l-1})}$. Then, uniformly in $n$, $z\in \partial D(0,r_{n, l})
$,
$w\in  \partial D(0,r_{n, l-1})$, and $B_{ n} \in \HH_{n, l}$,
\begin{equation}
\PPP^z (B_{ n} \,\big |\,X_{T_{  \partial D(0,r_{n, l-1})}}=w) = (1
+O(n^{ -3})) \PPP^{ z} (B_{ n}).
\label{m1.5euj}
\end{equation}  Furthermore, uniformly in $n$, $z,z'\in \partial D(0,r_{n, l})$,
and sequences 
$B_{ n}\in\mathcal{U}(\{ \HH_{n, l}\})$
\begin{equation}
\PPP^{ z} (B_{ n}) = (1 +O( n^{ -3})) \PPP^{ z'} (B_{n})
\,.
\label{m1.5evj}
\end{equation}
\end{lemma}

When we say that a statement such as (\ref{m1.5evj}) is uniform in  sequences 
$B_{ n}\in\mathcal{U}(\{ \HH_{n, l}\})$, we mean that there exists a $0<d<\ff$ 
such that if uniformly in $x,x'\in \partial D(0,r_{n, l+1})$
\begin{equation}
(1 - cn^{ -4}) \PPP^{ x'} (B_{ n})\leq \PPP^{ x} (B_{ n}) \leq (1 +cn^{ -4}) \PPP^{ x'}
(B_{ n})
\,.
\label{m1.5unifex}
\end{equation} then uniformly in $z,z'\in \partial D(0,r_{n, l})$
\begin{equation}
 (1 - dcn^{ -3}) \PPP^{ z'} (B_{
n})\leq\PPP^{ z} (B_{ n}) \leq  (1 + dcn^{ -3}) \PPP^{ z'} (B_{
n})
\,.
\label{m1.5evjf}
\end{equation}

\noindent{\bf Proof of Lemma \ref{lemprobends}:} Fixing $z \in  \partial D(0,r_{
n,l})$ it suffices to consider
$B_{ n} \in \HH_{ n,l}$ for which $\PPP^z(B_{ n})>0$. Fix such a set $B_{ n}$ and
a point
$w\in
\partial D(0,r_{ n,l-1})$. Using the notation introduced right before the
statement of our Lemma, for any $i
\geq 1$, we can write
\[\{ B_{ n}, I=i \}=\{ B_{n, i},
\tau_{2i}<\bar{\tau}\}\cap(\{I=0\}\circ\th_{\tau_{2i}})\]
 for some $B_{n, i}\in\FF_{ \tau_{2i}}=\{ A\,|\,A\cap \}$, 
 so by the strong Markov property at $\tau_{2i}$,
$$
P^z [ X_{\bar\tau}=w ; B_{ n}, I=i ] =
\E^z \left[ P^{X_{\tau_{2i}}} (X_{\bar\tau}=w,\, I=0) ; B_{n, i}, \tau_{2i} <
{\bar\tau} \right] \,,
$$ and
$$
\PPP^z \(B_{ n}, I=i \) = \E^z \left[ 
P^{X_{\tau_{2i}}} (I=0)\, ; B_{n, i}, \tau_{2i} < {\bar\tau} \right]
   \;.
$$ Consequently, for all $i \geq 1$, 
\begin{eqnarray} &&
P^z [ X_{\bar\tau}=w;  B_{ n}, I=i ]\label{ofer-star2}\\ &&
 \geq \PPP^z \(B_{ n}, I=i \)
\inf_{x\in \partial D(0,r_{n, l})}
\frac{P^{x} \( X_{\bar\tau}=w;\, I=0 \)} {P^x \(I=0\)}
\,. \nonumber
\end{eqnarray} 
Necessarily $\PPP^z(B_{ n}\, |I\,=0)\in \{0,1\}$
and is independent of $z$ for any $B_{ n} \in \HH_{n, l}$, implying that
\req{ofer-star2} applies for $i=0$ as well. By 
\req{p1.5h}, \req{p1.5} and \req{p1.2bfa}  there exists $c<\infty$
such that for any $z,x\in \partial D(0,r_{ n,l})$ and $w \in   \partial D(0,r_{
n,l-1})$,
\[
\frac{P^{x} \( X_{\bar\tau}=w;\, I=0 \)} {\E^x \(I=0\)}
 \geq (1-c n^{ -3})  H_{D(0,r_{n, l-1})^{ c}}(z,w)\,.
\]
 Hence, summing \req{ofer-star2} over $I=0,1,\ldots$, we get that
\[
   P^z \left[ X_{\bar\tau}=w, B_{ n} \right] \geq (1 -c n^{
-3})\PPP^z (B_{ n}) H_{ D(0,r_{n, l-1})^{ c}}(z,w)\,.
\] A similar argument shows that
$$
P^z \left[  X_{\bar\tau}=w, B_{ n} \right]
\leq (1+c n^{ -3}) \PPP^z \(B_{ n} \) H_{ D(0,r_{n, l-1})^{ c}}(z,w)\,,
$$ and we thus obtain (\ref{m1.5euj}).

By the Markov property at $\tau_1$, for any $z\in   \partial D(0,r_{n, l})$,
\bas &&
\PPP^{ z} (B_{ n})= \PPP^{z}(B_{ n}, I=0)
\\ &&\hspace{ 1in}+\sum_{ x\in  \partial D(0,r_{n, l+1})}H_{ D(0,r_{
n,l+1})\cup   D(0,r_{n, l-1})^{ c}}(z,x)\PPP^{x}\(B_{ n}\).
\eas 
The term involving $\{B_{ n} , I=0\}$ is dealt with
by (\ref{ps.3}) and by (\ref{m1.5unif}). For any fixed $x_{ 0}\in \partial D(0,r_{n,
l+1})$
\bea
&&
\sum_{ x\in  \partial D(0,r_{n, l+1})}H_{ D(0,r_{
n,l+1})\cup   D(0,r_{n, l-1})^{ c}}(z,x)\PPP^{x}\(B_{ n}\)\label{m1.5unifclar}\\
&&
=(1 +O( n^{ -4})) \PPP^{x_{ 0}}\(B_{ n}\)\sum_{ x\in  \partial
D(0,r_{n, l+1})}H_{ D(0,r_{ n,l+1})\cup   D(0,r_{n, l-1})^{ c}}(z,x)\nn\\
&&
=(1 +O( n^{ -4})) \PPP^{x_{ 0}}\(B_{ n}\)P^{ z}( T_{D(0,r_{ n,l+1}) }<
T_{D(0,r_{ n,l-1})^{ c} })\nn
\eea
so that (\ref{m1.5evj}) follows by  (\ref{ps.5}).
\qed 

Building upon Lemma \ref{lemprobends} we quantify the
 independence between the
$\si$-algebra $\GG_{n, l}^x$ of excursions from
$\partial D(x,r_{n, l-1})$ to $\partial D(x,r_{n, l})$ introduced in the
previous section and the
$\sigma$-algebra
$\HH_{n, l}^x (m)$ of excursions from $\partial D(x,r_{n,  l+1})$ to
$\partial D(x,r_{ n,l})$ during  the first $m$ excursions from
 $\partial D(x,r_{ n,l})$ to $\partial D(x,r_{ n,l-1})$. To this end, fix  
$x \in \Z^2$, let $\taup_0=0$ and for $i=1,2,\ldots$ recall that
\begin{eqnarray*}
\tau_{i} & = & \inf \{ t \geq \taup_{i-1}  :\; X_t \in \partial D(x,r_{n, l}) \} \,, \\
\taup_{i} & = & \inf \{ t \ge \tau_{i} :\; X_t \in  \partial D(x,r_{n, l-1})\}
\end{eqnarray*}  and that $\GG_{n, l}^x$ is the $\sigma$-algebra generated by the
excursions
$\{ e^{(j)}, j =1,\ldots \}$, where
$e^{(j)}=\{ X_{t} : \taup_{j-1}\leq t \leq \tau_{j} \}$ is the $j$-th excursion from
$\partial D(x,r_{ n,l-1})$ to $\partial D(x,r_{n, l})$ (so for $j=1$ we do begin at
$t=0$). 

We denote by
$\HH_{ n,l}^x (m)$ the $\sigma$-algebra generated by all excursions from
$\partial D(x,r_{ n, l+1})$ to
$\partial D(x,r_{ n,l})$ from time $\tau_1$ until time $\taup_{m}$. In more detail,
for each $j=1,2,\ldots,m$ let $\sip_{j,0}=\tau_{j}$ and for $i=1,\ldots$ define
\begin{eqnarray*}
\sii_{j,i} & = & \inf \{ t \geq \sip_{j,i-1}  :\; X_t \in
\partial  D(x,r_{ n, l+1})  \}\,, \\
\sip_{j,i} & = &
   \inf \{ t \ge \sii_{j,i} :\; X_t \in  \partial D(x,r_{n, l})\}\,.
\end{eqnarray*} Let $v_{j,i}=\{ X_{t} : \sii_{j,i}\leq t \leq \sip_{j,i} \}$ and
$Z^{j}=\sup \{ i \geq 0 \,:\,  \sip_{j,i}<\taup_{j}\}$. Then, $\HH_{n, l}^x(m)$ is the
product $\sigma$-algebra generated by the $\si$-algebras $\HH_{n, l,j}^x =
\sigma(v_{j,i}, i=1,\ldots,Z^{j})$ of the excursions between times $\tau_j$ and
$\taup_j$, for $j=1,\ldots,m$. 

We denote by $\mathcal{U}(\{ \HH_{ n,l}^x (m)\} )$ the collection of  sequences of
sets
$H_{ n}\in\HH_{n, l}^x (m)$ of the form $H_{ n}=H_{ n,1} \cap H_{ n,2} \cap
\cdots
\cap H_{ n,m}$, with $H_{ n,j} \in \mathcal{U}(\{ \HH^x_{n,l,j}\})$ for
$j=1,\ldots,m$.

\begin{lemma}\label{cond-ind} Uniformly over all $l\leq n$,
$m \leq ( n\log n)^{2}$,
$x,y_0,y_1 \in \Z^2$ and $H_{ n} \in \mathcal{U}(\{ \HH_{ n,l}^x (m)\})$,
\begin{eqnarray} && (1-O( m n^{ -3}) )  \PPP^{y_1} (H_{ n})
\leq
\PPP^{y_0} (H_{ n} \, |\, \GG_{n, l}^x
)\label{new1.2e}\\ &&
\hspace{ 2in}\leq (1+O( m n^{ -3}) )  \PPP^{y_1} (H_{ n})
  \,. \nonumber
\end{eqnarray}
\end{lemma}

\noindent{\bf Proof of Lemma \ref{cond-ind}:} We can write 
$H_{ n}=H_{ n,1} \cap H_{ n,2} \cap
\cdots
\cap H_{ n,m}$, with $H_{ n,j} \in \mathcal{U}(\{ \HH^x_{n,l,j}\})$ for
$j=1,\ldots,m$.  Conditioned upon
$\GG_{n, l}^x$  the events $H_{ n,j} $ are independent. Further, each $H_{ n,j} $
then has the conditional law of an event $B_{ n,j} $ in the collection $
\mathcal{U}(\{ \HH_{n,l}\})$ of Lemma
\ref{lemprobends}, for some random 
$z_j =X_{\tau_j}-x \in \partial D(0,r_{n, l})$
 and
$w_j =X_{\taup_j}-x \in \partial D(0,r_{n, l-1})$, both measurable on $\GG_{
n,l}^x$. By our conditions, the uniform estimates (\ref{m1.5euj}) and
(\ref{m1.5evj}) yield that for any fixed $z'\in \partial D(0,r_{ n,l})_{ s}$,
\beqn{ofer-star1} 
\PPP^{y_0} (H_{ n} \, |\, \GG_{n, l}^x )&&=\PPP^{y_0} (
\cap _{ j=1}^{ m}(H_{ n,j}) \, |\,
\GG_{n, l}^x )\\ && = \prod_{j=1}^m \PPP^{z_j}(B_{ n,j}
\, | \, X_{T_{\partial D(0,r_{n, l})^{ c}}} =w_j ) \nn \\ && = \prod_{j=1}^{m} \,(1 +
O( n^{ -3}))\PPP^{z_j}(B_{ n,j})\nn \\ && = (1 +O(n^{ -3})
)^{m}
\prod_{j=1}^{m} \,\PPP^{z'}(B_{ n,j} )\,.\nn
\eeqn Since $m \leq ( n\log n)^{2}$ and the right-hand side of \req{ofer-star1}
neither depends on
$y_0 \in \Z^2$ nor on the extra information in $\GG_{n, l}^x$, we get
\req{new1.2e} by averaging over $\GG_{n, l}^x$.
\qed

\noindent{\bf Proof of Lemma \ref{cond-cor}:} For $j=1,2,\ldots$ and
$i=l+1,l+2,\ldots,n$, let $Z_i^j$ denote the number of excursions from $\partial
D(y,r_{n,i-1})$ to $\partial D(y,r_{ n,i})$ by the random walk during the time
interval
$[\tau_j,\taup_j]$. Using (\ref{44.8}), the event
$$ H_{ n}=\{ \sum_{j=1}^{m_l} Z_i^j = m_i : i=l+1,l+2,\ldots,n \}
$$ can be written as a disjoint union of events in the collection 
${\mathcal U}(\{ \HH_{n, l}^y (m_l)\})$ of Lemma
\ref{cond-ind}. It is easy to verify that starting at any $x \notin D(y,r_{n,l})$, when
the event
$\{N^y_{n,l}=m_l\} \in \GG^y_{n,l}$ occurs, it implies that
$N^y_{n,i}=\sum_{j=1}^{m_l} Z_i^j$ for $i=l+1,l+2,\ldots,n$. Thus,
\begin{eqnarray} 
\PPP^{x_0} ( \Gamma_{ n,l}^y \, |{\mathcal
G}_{n,l}^y) {\bf 1}_{\{N^y_{n,l}=m_l\}}=
\PPP^{x_0} ( H_{ n}\, |{\mathcal G}_{n,l}^y) {\bf
1}_{\{N^y_{n,l}=m_l\}}\,. \label{hole0}
\end{eqnarray}  With $m_l/(n^2 \log n)$ bounded above, by \req{new1.2e} we
have, uniformly in 
$y \in \Z^2$ and $x_0,x_1 \in \Z^2 \setminus D(y,r_{n,l})$,
\begin{equation}
\label{hole1}
\PPP^{x_0} ( H_{ n}\, |{\mathcal G}_{n,l}^y)  =
(1+O( n^{-1}
 \log n))  
\PPP^{x_1} ( H_{ n})\,.
\end{equation}  Hence,
\begin{equation} 
\PPP^{x_0} ( \Gamma_{ n,l}^y\, |{\mathcal
G}_{n,l}^y)  {\bf 1}_{\{N^y_{n,l}=m_l\}} =(1+O( n^{-1}  \log
n))  
\PPP^{x_1} ( H_{ n}) {\bf 1}_{\{N^y_{n,l}=m_l\}}\,.
\label{hole2a}
\end{equation}
 Taking $x_0=x_1$ and averaging, one has
\begin{equation}
\PPP^{x_1} ( \Gamma_{ n,l}^y \, |
N^y_{n,l}=m_l)= (1+O(n^{-1}  \log n))  
\PPP^{x_1} ( H_{ n})\label{hole2}
\end{equation}
Hence, 
\begin{eqnarray}&& 
\PPP^{x_1} ( \Gamma_{ n,l}^y \, |
N^y_{n,l}=m_l){\bf 1}_{\{N^y_{n,l}=m_l\}}\label{hole2f}\\&&=(1+O(n^{-1} 
\log n))  
\PPP^{x_1} ( H_{ n}){\bf 1}_{\{N^y_{n,l}=m_l\}}
\nn \\&&=(1+O(n^{-1} 
\log n)) \PPP^{x_0} ( \Gamma_{ n,l}^y\, |{\mathcal
G}_{n,l}^y)  {\bf 1}_{\{N^y_{n,l}=m_l\}} \nonumber
\end{eqnarray}
where we used (\ref{hole2a}) for the last equality.
 Using that
${\{N^y_{n,l}=m_l\}}\subset {\mathcal G}_{n,l}^y\,$, this
is \req{new1.3e}.
\qed
\vspace{2mm}

\noindent{\bf Acknowledgements}
I am grateful to Endre Cs\'{a}ki and  Ant\'{o}nia F\"{o}ldes for many helpful
comments.

\bigskip
\noindent
\begin{tabular}{lll} 
  & Jay Rosen & \\
  & Department of Mathematics&\\
 &College of Staten Island, CUNY&\\
 &Staten Island, NY 10314&\\ &jrosen3@earthlink.net&
\end{tabular}

\end{document}